\documentclass{article}
\usepackage[utf8]{inputenc}

\usepackage{amsmath}
\usepackage{amssymb}
\usepackage{mathtools}
\usepackage{graphicx}

\title{Gaussian Integers, Rings, Finite Fields,\\ and the Magic Square of Squares}
\author{O. M. Cain \\ onnomawcain@gmail.com}
\date{April 2019}

\begin{document}

\maketitle

\textbf{Abstract.} We show the $3\times3$ magic square of squares problem equivalent to solving quartic polynomials with certain factorization constraints over an abelian extension of the rationals. We analyze a particular case in which said extension is assumed to be the Gaussian integers resulting a new search method. Additionally, the magic square of squares is analyzed over finite fields and rings of the form $\mathbb{Z}/n\mathbb{Z}$ resulting in some conjectures enumerating the rings and finite fields in which a magic square of squares can be constructed. Code is made available.

\section{Background}

The construction of a $3\times3$ magic square of squares -- sometimes called simply \textit{the magic square of squares problem} -- is defined to be $9$ distinct squared integers placed in a $3\times3$ grid, 
$$\begin{bmatrix}
a^2&b^2&c^2\\
d^2&e^2&f^2\\
g^2&h^2&i^2\\
\end{bmatrix},$$
such that the sums of the elements in each row, column, and the two main diagonals sum to the same total. That is to say for some integer total $T$ we have
$$a^2+b^2+c^2=d^2+e^2+f^2=g^2+h^2+i^2=T,$$
$$a^2+d^2+g^2=b^2+e^2+h^2=c^2+f^2+i^2=T,$$
$$\text{and}\quad a^2+e^2+i^2=g^2+e^2+c^2=T.$$
In total there are eight sums to be satisfied. A lot has been written up on this fabled object [1] but it is currently unknown if any solution exists. Some 'near misses' have been found such as the \textit{Parker Square} [2],
$$\begin{bmatrix}
29^2&1^2&47^2\\
41^2&37^2&1^2\\
23^2&41^2&29^2\\
\end{bmatrix},$$
in which seven of the eight sums add to $3051$ but the eighth sums to $4107\not= 3051$. The Parker Square also unfortunately has two duplicate entries.\\
\indent The problem has been shown related to congruent numbers, elliptic curves, and right triangles with rational sides [3] (which isn't terribly surprising since each of those are tightly related to the others).

\section{The Magic Hourglass of Squares}

At simpler -- and still unsolved -- problem is whether there exists a magic hourglass of squares. That is, $7$ squared integers, 
$$\begin{bmatrix}
a^2&b^2&c^2\\
-&e^2&-\\
g^2&h^2&i^2\\
\end{bmatrix},$$
satisfying
$$a^2+e^2+i^2=b^2+e^2+h^2=c^2+e^2+g^2=T$$
$$\text{and}\quad a^2+b^2+c^2=g^2+h^2+i^2=T$$
for some integer total $T$. There are two overlapping magic hourglasses of squares in any magic square of squares. We will warm up by first observing \\\\
\textbf{Theorem 2.1:} The total of any magic hourglass of squares $T$ is $3$ times the central entry $e^2$.\\\\
\textbf{Proof:} Observe by rearranging
$$3T = (a^2+e^2+i^2)+(b^2+e^2+h^2)+(c^2+e^2+g^2)$$
$$=(a^2+b^2+c^2)+(e^2+e^2+e^2)+(g^2+h^2+i^2)=2T+3e^2.$$
The Theorem follows immediately man.\\
\indent \hfill$\square$\\
\textbf{Corollary 2.1:} For any magic hourglass of squares, entries as given above, we have
$$a^2+i^2=b^2+h^2=c^2+g^2=2e^2.$$
\textbf{Proof:} From Theorem 2.1 we have $a^2+e^2+i^2=3e^2$ or equivalently $a^2+i^2=2e^2$. The same follows for the other two sums containing the central entry.\\
\indent \hfill$\square$\\
Solutions to $r^2+t^2=2s^2$ over the integers have been known for a long time [4]. Note this is equivalent to an arithmetic progression of $3$ squares since $r^2-s^2=s^2-t^2$. Specifically\\\\
\textbf{Theorem 2.2:} For every integer solution to $r^2+t^2=2s^2$ there exist 3 integer parameters $m,n,$ and $k$ such that 
$$r = k(m^2+2mn-n^2),\quad s=k(m^2+n^2),\quad t=k(m^2-2mn-n^2),$$
$$\quad r^2-s^2=s^2-t^2=4mn(m+n)(m-n)k^2.$$
\textbf{Proof:} Left undone.\\\\
The converse of Theorem 2.2 is also true: any such $m,n,$ and $k$ solve $r^2+t^2=2s^2$. \\
\indent The parameter $k$ in is only a scaling factor. It is $m$ and $n$ which are really ``doing the work". For example, we may pick $m=3$ and $n=2$ to obtain
$$17^2-13^2=13^2-7^2=120.$$
\indent An existing guess-and-check method of hourglass searching makes use of this fact [1]. The `guess' is made by picking a number expressible as the sum of two squares in at least 3 ways, say 
$$1105=33^2+4^2=32^2+9^2=24^2+23^2,$$
and making use of Theorem 2.2 to generate the 3 sums of the hourglass which include the central entry. In this case we get
$$367^2+1519^2=1337^2+809^2=1057^2+1151^2=2(1105^2)$$ and the corresponding hourglass
$$\begin{bmatrix}
367^2&1337^2&1151^2\\
-&1105^2&-\\
1057^2&809^2&1519^2\\
\end{bmatrix}$$
in which the 3 of the 5 sums have the same total. The top and bottom rows are different.

\section{Factorization Over the Gaussian Integers}

Theorem 2.2 has a nice reinterpretation as\\\\
\textbf{Lemma 3.1:} For every integer solution to $r^2+t^2=2s^2$ there exists a complex number $\omega\in\mathbb{Z}[i,\sqrt{k}]$ such that $k$ is an integer and
$$r=\text{Re}[\omega^2]+\text{Im}[\omega^2],\quad s=\omega\overline{\omega},\quad t=\text{Re}[\omega^2]-\text{Im}[\omega^2],$$
$$r^2-s^2=s^2-t^2=\text{Im}[\omega^4],\quad\text{and}\quad rt=\text{Re}[\omega^4].$$
\textbf{Proof:} Choose $\omega=\sqrt{k}(m+ni)$. All but the last equation follow from Theorem 2.2. Noting that $\text{Re}[xy]=\text{Re}[x]\text{Re}[y]-\text{Im}[x]\text{Im}[y]$, the last equation is derived easily.\\
\indent \hfill$\square$\\
\indent To simplify notation, we let the function $\chi$ take a complex number $\omega $ to $(r,s,t)$ as defined above. Meaning specifically $\chi(\omega)=(\text{Re}[\omega^2]+\text{Im}[\omega^2],\omega\overline{\omega},\text{Re}[\omega^2]-\text{Im}[\omega^2])$\\\\
\textbf{Theorem 3.2:} A magic hourglass of squares, 
$$\begin{bmatrix}
a^2&b^2&c^2\\
-&e^2&-\\
g^2&h^2&i^2\\
\end{bmatrix},$$ 
exists if and only if there exists $\alpha,\beta,\gamma\in K=\mathbb{Z}[i,\sqrt{A},\sqrt{B},\sqrt{C}]$ with $A,B,$ and $C$ integers such that
$$\alpha^4+\beta^4+\gamma^4\in\mathbb{Z},\quad\alpha^2,\beta^2,\gamma^2\in\mathbb{Z}[i],\quad \alpha\overline{\alpha}=\beta\overline{\beta}=\gamma\overline{\gamma},$$
$$(a,e,i)=\chi(\alpha),\quad (b,e,h)=\chi(\beta),\quad (c,e,g)=\chi(\gamma),$$
with $\alpha^4,\beta^4,$ and $\gamma^4$ distinct and strictly complex.\\\\
\textbf{Proof:} From Corollary 2.1 note that in any magic hourglass of squares $(a^2,e^2,i^2),$ $(b^2,e^2,h^2),$ and $(c^2,e^2,g^2)$ each form an arithmetic progression. Applying Lemma 3.1 to each we have $\alpha,\beta, $ and $\gamma$ parametrize the 3 progressions respectively -- in the same way that $\omega$ parametrized $r^2,s^2,$ and $t^2$. Note the equivalence of $\alpha\overline{\alpha},\beta\overline{\beta},$ and $\gamma\overline{\gamma}$ follows from the fact that the 3 progressions have the same middle element. \\
\indent By Theorem 2.2 and Lemma 3.1 we see that the sums $a^2+e^2+i^2,b^2+e^2+h^2,$ and $c^2+e^2+g^2$ are all equal to $3e^2$ if and only if we have the former parametrization by $\alpha,\beta,$ and $\gamma$. Next we must show that $a^2+b^2+c^2=3e^2$ (the logic for $g^2+h^2+i^2$ is identical). Rearranging, we have
$$(a^2-e^2)+(b^2-e^2)+(c^2-e^2)=0.$$
By Lemma 3.1 again, we see this is equivalent to $\text{Im}[\alpha^4]+\text{Im}[\beta^4]+\text{Im}[\gamma^4]=0.$ Which in turn is equivalent to
$$\alpha^4+\beta^4+\gamma^4\in\mathbb{Z}.$$
\indent The restriction $\alpha^4\not\in \mathbb{Z}$ ensures that $a^2,e^2,$ and $i^2$ are distinct (an analogous statement can be said for $\beta$ and $\gamma$). Restricting $\alpha^4,\beta^4,$ and $\gamma^4$ to be distinct ensures that the remaining elements of the hourglass are distinct as well.\\
\indent As for the reverse direction, it can be seen easily from Theorem 2.2 and Lemma 3.1 that any choice of $\alpha, \beta,\gamma$ within the restrictions yields a valid magic hourglass.\\
\indent \hfill$\square$\\
\textbf{Observation 3.1:} The Galois group of $K=\mathbb{Q}(i, \sqrt{A},\sqrt{B},\sqrt{C})$ is a product of cyclic groups of order $2$ and is therefore abelian. It follows (we think?) that K is a subring of the ring of integers of some cyclotomic field extension of the rationals. That sounds useful but we weren't able to do anything with it.\\\\
\textbf{Observation 3.2:} For a full magic square -- as opposed to a measly hourglass -- we obtain two inclusions
$$\alpha^4+\beta^4+\gamma^4,\alpha^4+\delta^4-\gamma^4\in\mathbb{Z}$$
where $\chi(\delta)=(d,e,f)$ and obeys a similar set of restrictions.

\section{Case $K=\mathbb{Z}[i]$}

The Gaussian integers $\mathbb{Z}[i]$ have unique factorization (see Theorem 6.6 in [6]). Thus if we take $\alpha,\beta,\gamma\in\mathbb{Z}[i]$ (meaning that the parameter $k$ as it appears in Lemma 3.1 is $1$ for each of $\alpha, \beta,$ and $\gamma$) then by 
\begin{equation}
\label{norms}
    \alpha\overline{\alpha}=\beta\overline{\beta}=\gamma\overline{\gamma}
\end{equation}
from Theorem 3.2, we see that $\alpha,\beta,$ and $\gamma$ have the same prime factorization up to conjugation of factors. This allows for some interesting formulations such as\\\\
\textbf{Theorem 4.1:} If there exists $x,y,z\in\mathbb{Z}[i]$ such that
$$\text{Im}[x^4y^4z^4]=-4\text{Im}[x^4]\text{Im}[y^4]\text{Im}[z^4]$$
and $x^4,y^4,$ and $z^4$ are both strictly complex and not real multiples of each other then there exists a magic hourglass of squares. \\\\
\textbf{Proof:} Choose
$$\alpha=\overline{x}yz,\quad \beta=x\overline{y}z,\quad\text{and}\quad \gamma=xy\overline{z}.$$
It follows immediately that $\alpha\overline{\alpha}=\beta\overline{\beta}=\gamma\overline{\gamma}$ and trivially that $\alpha^2,\beta^2,\gamma^2\in\mathbb{Z}[i]$. To see that $\alpha^4+\beta^4+\gamma^4$ is real we must first note the identity
$$\text{Im}[XYZ]=$$
$$\text{Re}[X]\text{Re}[Y]\text{Im}[Z]+\text{Re}[X]\text{Im}[Y]\text{Re}[Z]+\text{Im}[X]\text{Re}[Y]\text{Re}[Z]-\text{Im}[X]\text{Im}[Y]\text{Im}[Z].$$
From this one derives
$$\text{Im}[x^4y^4z^4]+4\text{Im}[x^4]\text{Im}[y^4]\text{Im}[z^4]=\text{Im}[\overline x^4y^4z^4]+\text{Im}[x^4\overline y^4z^4]+\text{Im}[x^4y^4\overline z^4]$$
$$=\text{Im}[\alpha^4]+\text{Im}[\beta^4]+\text{Im}[\gamma^4]$$
to see that $\text{Im}[\alpha^4+\beta^4+\gamma^4]=0$. \\
\indent To see that none of $\alpha^4,\beta^4,$ or $\gamma^4$ is real, suppose that one of them, say $\alpha^4$, is. Then necessarily $\text{Im}[\beta^4]=-\text{Im}[\gamma^4]$. Along with $\beta\overline{\beta}=\gamma\overline{\gamma}$, this implies that $\gamma^4=\overline{\beta}^4$ or equivalently $\gamma^2=\pm\overline{\beta}^2$. But if so then $\beta^2\gamma^2=\pm x^4(y\overline{y})^2(z\overline{z})^2$ is real implying in turn that $x^4$ is real; a contradiction.\\
\indent Lastly, to see that $\alpha^4,\beta^4,$ and $\gamma^4$ are distinct suppose that two are not, say $\alpha^4=\beta^4$. It follows that $\alpha^2=\pm\beta^2$ and therefore that $\overline{x}^2y^2=\pm x^2\overline{y^2}$. Scaling both sides by $x^2y^2$ we obtain $|x|^4y^4=\pm x^4|y|^4$ which implies that $x^4$ and $y^4$ are real multiples of each other; a contradiction.\\
\indent \hfill$\square$\\
It's not clear whether the converse of Theorem 4.1 is true. But a slightly weaker statement can be obtained replacing $x^4,y^4,$ and $z^4$ with $x^2,y^2,$ and $z^2$ respectively.\\\\
\textbf{Theorem 4.2:} If there exists a magic hourglass of squares, then there exists $x,y,z\in\mathbb{Z}$ such that 
$$\text{Im}[x^2y^2z^2]=-4\text{Im}[x^2]\text{Im}[y^2]\text{Im}[z^2]$$
and $x^2,y^2,$ and $z^2$ are both strictly real and not real multiples of each other.\\\\
\textbf{Proof:} Left undone.\\\\
\textbf{Observation 4.1:} Theorem 4.1 can be used to search for magic hourglasses. It's known that for any Gaussian integer $\omega=m+ni$ the value $\text{Im}[\omega^4]=4mn(m+n)(m-n)$ is always divisible by $24$ [4]. Thus $4\text{Im}[x^4]\text{Im}[y^4]\text{Im}[z^4]$ is divisible by $4\cdot 24^3$. One may then iterate over integers $A$ and $B$ assigning $xyz=2^7A+3^2B$ (which will necessarily satisfy the former divisibility constraint) and then check for each possible factorization, $x,y,$ and $z$, whether the constraint of Theorem 4.1 is also satisfied. It's not clear if there's a systematic way to choose integers $A$ and $B$ so that $xyz$ is easy to factor or -- perhaps alternatively -- has many factorizations.

\section{The Smallest Non-Parker Finite Field}

At this point we move on from the traditional magic square of squares problem. Specifically, we will attempt to find solutions over elements of finite fields instead of over the integers. We call a finite field -- or ring in general -- \textit{Parker} if no magic square of square can be formed; that is, if no nine distinct squared elements can be found satisfying the constraints given in the Section 1. If a traditional magic square of squares does not exist, then $\mathbb{Z}$ is Parker. It is hoped that by enumerating which finite fields and rings are Parker, it may become clear whether $\mathbb{Z}$ is Parker.\\
\indent The field $\mathbb{F}_{29}$, for example, is not Parker since
\begin{equation*}
\begin{tabular}{c|c|c}
  $9^2$ & $11^2$ & $1^2$ \\      \hline
  $6^2$ & $0^2$ & $14^2$ \\      \hline
  $12^2$ & $16^2$ & $8^2$
\end{tabular}\ 
\end{equation*}
is a valid magic square of squares evaluated  mod $29$.\\
\indent As a guiding inquiry, we ask what is the non-Parker finite field of smallest order? We have seen already that $\mathbb{F}_{29}$ is non-Parker, so there are only finitely many fields left to check. Remember that for any prime, $p$, or prime power, $p^r$, there exists a finite field -- isomorphically unique -- having that order (see Theorem 15.7.3 in [7]).\\ 
\indent The finite fields with orders less than $29$ are 
$$\mathbb{F}_{2}, \mathbb{F}_{3}, \mathbb{F}_{4}, \mathbb{F}_{5}, \mathbb{F}_{7}, \mathbb{F}_{8}, \mathbb{F}_{9}, \mathbb{F}_{11}, \mathbb{F}_{13}, \mathbb{F}_{16}, \mathbb{F}_{17}, \mathbb{F}_{19}, \mathbb{F}_{23}, \mathbb{F}_{25}, \ \text{and}\ \mathbb{F}_{27}. $$
One could quickly check which of these, if any, are Parker with symbolic computation. We will however proceed by whittling down this list, lemma by lemma, in the hope of also getting a better sense of the structure of magic squares of squares. If any reader would like a more rigorous treatement of the cases $\mathbb{F}_2^r$ and $\mathbb{F}_3^r$ than what we can offer, they are referred to [10]. We start with \\\\
\textbf{Lemma 5.1:} All $3\times3$ magic squares over finite fields of even order have duplicate entries.\\\\
\textbf{Proof:} All finite fields can be formed by taking a quotient of a polynomial ring over a finite field of prime order (see Theorem 15.7.3 of [2]). For example, 
$$\mathbb{F}_{9}\cong \mathbb{F}_3[x]/(x^2+1).$$
This means every element of a field of even order is a sum of distinct powers of $\overline{x}$. For example, the entries of $\mathbb{F}_8\cong \mathbb{F}_2[x]/(x^3+x+1)$ are 
$$0,\ \ 1,\ \ \overline{x},\ \ \overline{x}+1,\ \ \overline{x}^2,\ \ \overline{x}^2+1,\ \ \overline{x}^2+\overline{x},\ \text{and}\ \overline{x}^2+\overline{x}+1.$$
Thus any magic square in $\mathbb{F}_{2^k}$ can be regarded as a linear combination of magic squares (added element-wise) over $\mathbb{F}_2$ where each such square is scaled by a power of $\overline{x}$. A magic square over $\mathbb{F}_{16}$,
\begin{equation*}
\text{take}\quad
\begin{tabular}{c|c|c}
  $1+\overline{x}+\overline{x}^3$ & $\overline{x}+\overline{x}^2$ & $1+\overline{x}^2+\overline{x}^3$ \\      \hline
  $\overline{x}+\overline{x}^2$ & $0$ & $\overline{x}+\overline{x}^2$ \\      \hline
  $1+\overline{x}^2+\overline{x}^3$ & $\overline{x}+\overline{x}^2$ & $1+\overline{x}+\overline{x}^3$
\end{tabular}\quad \text{for example,}
\end{equation*}
can be decomposed as an element-wise sum of $4$ magic squares over $\mathbb{F}_2$,
\begin{equation*}
\begin{tabular}{c|c|c}
  $1$ & $0$ & $1$ \\      \hline
  $0$ & $0$ & $0$ \\      \hline
  $1$ & $0$ & $1$
\end{tabular}\ +\ \overline{x}\cdot 
\begin{tabular}{c|c|c}
  $1$ & $1$ & $0$ \\      \hline
  $1$ & $0$ & $1$ \\      \hline
  $0$ & $1$ & $1$
\end{tabular}\ +\ \overline{x}^2\cdot 
\begin{tabular}{c|c|c}
  $0$ & $1$ & $1$ \\      \hline
  $1$ & $0$ & $1$ \\      \hline
  $1$ & $1$ & $0$
\end{tabular}\ +\ \overline{x}^3\cdot 
\begin{tabular}{c|c|c}
  $1$ & $0$ & $1$ \\      \hline
  $0$ & $0$ & $0$ \\      \hline
  $1$ & $0$ & $1$
\end{tabular}\ .
\end{equation*}
\indent Naturally, we proceed by analyzing magic squares in $\mathbb{F}_2$. There are $8$:
\begin{equation*}
\begin{tabular}{c|c|c}
  $0$ & $0$ & $0$ \\      \hline
  $0$ & $0$ & $0$ \\      \hline
  $0$ & $0$ & $0$
\end{tabular}\ ,\ 
\begin{tabular}{c|c|c}
  $1$ & $1$ & $0$ \\      \hline
  $1$ & $0$ & $1$ \\      \hline
  $0$ & $1$ & $1$
\end{tabular}\ ,\ 
\begin{tabular}{c|c|c}
  $0$ & $1$ & $1$ \\      \hline
  $1$ & $0$ & $1$ \\      \hline
  $1$ & $1$ & $0$
\end{tabular}\ ,\ 
\begin{tabular}{c|c|c}
  $1$ & $0$ & $1$ \\      \hline
  $0$ & $0$ & $0$ \\      \hline
  $1$ & $0$ & $1$
\end{tabular}\ ,\ 
\end{equation*}
\begin{equation*}
\begin{tabular}{c|c|c}
  $1$ & $1$ & $1$ \\      \hline
  $1$ & $1$ & $1$ \\      \hline
  $1$ & $1$ & $1$
\end{tabular}\ ,\ 
\begin{tabular}{c|c|c}
  $0$ & $0$ & $1$ \\      \hline
  $0$ & $1$ & $0$ \\      \hline
  $1$ & $0$ & $0$
\end{tabular}\ ,\ 
\begin{tabular}{c|c|c}
  $1$ & $0$ & $0$ \\      \hline
  $0$ & $1$ & $0$ \\      \hline
  $0$ & $0$ & $1$
\end{tabular}\ ,\ \text{and}\ \ 
\begin{tabular}{c|c|c}
  $0$ & $1$ & $0$ \\      \hline
  $1$ & $1$ & $1$ \\      \hline
  $0$ & $1$ & $0$
\end{tabular}\ .
\end{equation*}
\indent To see that these are the only possible $8$, we resort to a previously published result [8] that any $3\times3$ magic square over the integers can be parametrized by $3$ integers, $A,B,$ and $C$, in the form
\begin{equation*}
\begin{tabular}{c|c|c}
  $A$ & $-A$ & $A$ \\      \hline
  $-A$ & $0$ & $A$ \\      \hline
  $0$ & $A$ & $-A$
\end{tabular}\ +\ 
\begin{tabular}{c|c|c}
  $0$ & $-B$ & $B$ \\      \hline
  $B$ & $0$ & $-B$ \\      \hline
  $-B$ & $B$ & $0$
\end{tabular}\ +\ 
\begin{tabular}{c|c|c}
  $C$ & $C$ & $C$ \\      \hline
  $C$ & $C$ & $C$ \\      \hline
  $C$ & $C$ & $C$
\end{tabular}\ .
\end{equation*}
\indent And now we arrive at the important bit. All $8$ squares have the same entry in the edge-middles and in opposite corners. Thus any magic square over a finite field of even order will have at most $4$ distinct entries.\\
\indent \hfill $\square$\\\\
\textbf{Corollary 5.1:} All finite fields of even order are Parker.\\
\indent\textbf{Proof:} By Lemma 3.1 any $3\times3$ magic square over a finite field has duplicate entries.\\
\indent \hfill $\square$\\\\
\indent Okay. First lemma down. We've eliminated $\mathbb{F}_2,\mathbb{F}_4,\mathbb{F}_8,$ and $\mathbb{F}_{16}$ from our list of candidates leaving
$$\mathbb{F}_{3}, \mathbb{F}_{5}, \mathbb{F}_{7}, \mathbb{F}_{9}, \mathbb{F}_{11}, \mathbb{F}_{13}, \mathbb{F}_{17}, \mathbb{F}_{19}, \mathbb{F}_{23}, \mathbb{F}_{25}, \ \text{and}\ \mathbb{F}_{27}. $$
The smaller fields can be eliminated easily.\\\\
\textbf{Lemma 5.2:} A finite field of odd order, $q$, has $\frac{q+1}{2}$ squares.\\
\indent \textbf{Proof:} This follows easily from the fact that $\mathbb{F}_q^\times$ is cyclic (for proof of which, we cite Artin [7] again; Theorem 15.7.3).\\
\indent \hfill $\square$\\\\
\textbf{Corollary 5.2:} The fields $\mathbb{F}_{3},\mathbb{F}_{5},\mathbb{F}_{7},\mathbb{F}_{9},\mathbb{F}_{11},$ and $\mathbb{F}_{13}$ are Parker.\\
\indent \textbf{Proof:} By Lemma 5.2, each of the fields in question have fewer than $9$ distinct squares. Therefore no $3\times3$ magic square of distinct squares can be formed.\\
\indent \hfill $\square$\\\\
\indent The remaining candidates are $\mathbb{F}_{17}, \mathbb{F}_{19}, \mathbb{F}_{23}, \mathbb{F}_{25}, \ \text{and}\ \mathbb{F}_{27}. $ Tackling these remaining fields will require poking into their structure somewhat deeper (either that or just performing a computation, but we decided not to take that route).\\\\
\textbf{Lemma 5.3:} Any non-Parker finite field contains either $4$ distinct solutions to $x^2+y^2=0$ with $x,y\not=0$ or $4$ distinct solutions to $x^2+y^2=2$ with $x^2,y^2\not=2$.\\
\indent\textbf{Proof:} By definition, a non-Parker field contains at least one magic square of distinct squares. We use the standard variables
\begin{equation*}
\begin{tabular}{c|c|c}
  $a^2$ & $b^2$ & $c^2$ \\      \hline
  $d^2$ & $e^2$ & $f^2$ \\      \hline
  $g^2$ & $h^2$ & $i^2$
\end{tabular}\ .
\end{equation*}
If $e=0$ then there are $4$ distinct solutions to $x^2+y^2=0$. If $e$ is nonzero, we use the fact that all elements in a field have inverses (this is one of the facts that defines what exactly a ``field" is). We scale our hypothetical magic square by $(e^{-1})^2$ which produces another magic square of distinct squares. The new square has a central entry of $1$. Thus there are $4$ distinct solutions to $x^2+y^2=2(1^2)=2$.\\
\indent \hfill $\square$\\\\
\textbf{Corollary 5.3:} The fields $\mathbb{F}_{19},\mathbb{F}_{23},$ and $\mathbb{F}_{27}$ are Parker.\\
\indent\textbf{Proof:} Let's count the solutions $x^2+y^2=0,2$ in each of of the fields in question. This took us roughly 10 minutes per field to do by hand (and was, in fact, further verified with computation [9]). There are no solutions to $x^2+y^2=0$. The respective solutions to $x^2+y^2=2$ are
$$\mathbb{F}_{19}:\quad 2^2+6^2=4^2+9^2=2, $$
$$\mathbb{F}_{23}:\quad 0^2+5^2=3^2+4^2=6^2+9^2=2, $$
$$\text{and}\quad \mathbb{F}_{27}:\quad (\overline{x})^2+(\overline{x}^2+2\overline{x})^2=(\overline{x}^2+2)^2+(\overline{x}+2)^2=(\overline{x}+1)^2+(\overline{x}^2+\overline{x})^2=2 $$
where $\overline{x}^3=\overline{x}+2$. None of the fields have $4$ distinct solutions and are thus each Parker by Lemma 5.3.\\
\indent \hfill $\square$\\\\
\indent We are left with $\mathbb{F}_{17}\ \text{and}\ \mathbb{F}_{25}$
each having $4$ distinct solutions to $x^2+y^2=0$.\\
\indent But this doesn't yet ensure $\mathbb{F}_{17}$ and $\mathbb{F}_{25}$ are non-Parker. We need one last lemma.\\\\
\textbf{Lemma 5.4:} Magic squares of distinct squares over a finite field with a central entry of $0$ are parametrized (up to scaling) by solutions to $\alpha^2-\beta^2=\beta^2-\gamma^2=1$ (i.e. three consecutive squares) satisfying $\{\alpha,\beta,\gamma\}\cap \{0,1,-1\}=\emptyset$.\\\\
\textbf{Proof:} Suppose we have a magic square of squares over a finite field with a central entry of $0$, 
\begin{equation*}
\begin{tabular}{c|c|c}
  $a^2$ & $b^2$ & $c^2$ \\      \hline
  $d^2$ & $0$ & $-d^2$ \\      \hline
  $-c^2$ & $-b^2$ & $-a^2$
\end{tabular}\ .
\end{equation*}
Again, we use the fact that all entries have inverses and scale the square by $(c^{-1})^2$. With the right change of variables, we obtain
\begin{equation*}
\begin{tabular}{c|c|c}
  $\beta^2$ & $-\alpha^2$ & $1$ \\      \hline
  $-\gamma^2$ & $0$ & $\gamma^2$ \\      \hline
  $-1$ & $\alpha^2$ & $-\beta^2$
\end{tabular}\ .
\end{equation*}
This square is manifestly magic if and only if 
$$\alpha^2-\beta^2-1=\beta^2-\gamma^2-1=0.$$
The lemma follows immediately.\\
\indent \hfill $\square$\\\\
\textbf{Corollary 5.4:} The fields $\mathbb{F}_{17}$ and $\mathbb{F}_{25}$ are Parker.\\\\
\textbf{Proof:} Up to isomorphism, the squares of $\mathbb{F}_{17}$ and of $\mathbb{F}_{25}$ are 
$$\{0,1,2,4,8,9,13,15,16\}$$
$$\text{and}\quad \{0,1,2,3,4,\overline{x}+1,\overline{x}+3,2\overline{x}+1,2\overline{x}+1,3\overline{x}+1,3\overline{x}+1,4\overline{x}+1,4\overline{x}+1\}$$
respectively. They manifestly have no three consecutive squares (excluding $0$ and $\pm1$). Thus, by our Lemma 5.4, they are Parker.\\
\indent \hfill $\square$\\
\indent We've finished the whole list of candidates. The final result of this section can be stated.\\\\
\textbf{Theorem 5.1:} $\mathbb{F}_{29}$ is the non-Parker field of smallest order.\\
\indent\textbf{Proof:} All finite fields of smaller order are Parker by Corollaries 5.1, 5.2, 5.3, and 5.4. We see that $\mathbb{F}_{29}$ is non-Parker by the aforementioned construction.\\
\indent \hfill $\square$

\section{Search Algorithm for Finite Fields}

We did our best in the previous section to work by argumentation. In this section we'll switch to empirical results obtained via a computer algebra system [9]. It will first be addressed how to generate all magic squares over a given field. In pseudo-code:\\\\
\textbf{Algorithm 6.1:}\\
\# Input: A finite field, $\mathbb{F}_q$.\\
\# Output: Set of all tuples $(a^2,b^2,...,i^2)$ forming magic squares over $\mathbb{F}_q$\\
\# $\quad$ up to scaling.\\
\textbf{function} msos\_field($\mathbb{F}_q$):\\
\indent SQUARES $\leftarrow \{x^2:x^2\in \mathbb{F}_q\}$\\
\indent MSOS $\leftarrow \{\}$\\
\indent $e\leftarrow 0 \quad\quad\text{\# first case.}$\\
\indent \textbf{for} $\{a^2,i^2\}\subset$ SQUARES:\\
\indent\indent \textbf{if} $a^2+i^2\not=2e^2$: \textbf{continue}\\
\indent\indent $c^2,g^2\leftarrow 1,-1$\\
\indent\indent $B\leftarrow 3e^2-a^2-c^2$ \\
\indent\indent $D\leftarrow 3e^2-a^2-g^2$ \\
\indent\indent $F\leftarrow 3e^2-c^2-i^2$ \\
\indent\indent $H\leftarrow 3e^2-g^2-i^2$ \\
\indent\indent \textbf{if} $\{B,D,F,H\}\not\subset$ SQUARES: \textbf{continue}\\
\indent\indent \textbf{if} $|\{a^2,B,c^2,D,e^2,F,g^2,H,i^2\}|\not=9$: \textbf{continue}\\
\indent\indent MSOS $\leftarrow \text{MSOS}\cup \{(a^2,B,c^2,D,e^2,F,g^2,H,i^2)\}$\\
\indent $e\leftarrow 1 \quad\quad\text{\# second case.}$\\
\indent SEQUENCES $\leftarrow \{\}$\\
\indent \textbf{for} $\{a^2,i^2\}\subset$ SQUARES:\\
\indent\indent \textbf{if} $a^2+i^2\not=2e^2$: \textbf{continue}\\
\indent\indent \textbf{for} $\{c^2,g^2\}\subset$ SEQUENCES:\\
\indent\indent\indent $B\leftarrow 3e^2-a^2-c^2$ \\
\indent\indent\indent $D\leftarrow 3e^2-a^2-g^2$ \\
\indent\indent\indent $F\leftarrow 3e^2-c^2-i^2$ \\
\indent\indent\indent $H\leftarrow 3e^2-g^2-i^2$ \\
\indent\indent\indent \textbf{if} $\{B,D,F,H\}\not\subset$ SQUARES: \textbf{continue}\\
\indent\indent\indent \textbf{if} $|\{a^2,B,c^2,D,e^2,F,g^2,H,i^2\}|\not=9$: \textbf{continue}\\
\indent\indent\indent MSOS $\leftarrow \text{MSOS}\cup \{(a^2,B,c^2,D,e^2,F,g^2,H,i^2)\}$\\
\indent\indent SEQUENCES $\leftarrow\text{SEQUENCES}\cup\{\{a^2,i^2\}\}$\\
\indent \textbf{return} MSOS\\\\
\indent Implementation in the SageMath language is available at [11].\\
\indent The correctness of the code follows roughly from some observations made while proving Lemmas 5.3 and 5.4. We neglect to give a rigorous proof, but will content ourselves noting the following intuition. \\
\indent Any magic square of squares over a finite field may be scaled into one of two forms according to whether the central entry is zero:
\begin{equation*}
\begin{tabular}{c|c|c}
  $a^2$ & $b^2$ & $1$ \\      \hline
  $d^2$ & $0$ & $f^2$ \\      \hline
  $-1$ & $h^2$ & $i^2$
\end{tabular}\quad\text{or}\quad
\begin{tabular}{c|c|c}
  $a^2$ & $b^2$ & $c^2$ \\      \hline
  $d^2$ & $1$ & $f^2$ \\      \hline
  $g^2$ & $h^2$ & $i^2$
\end{tabular}\ .
\end{equation*}
Observe that in the first case, the square is determined by just $a^2$ and, in the second case, by $a^2$ and $c^2$. Both cases are labeled accordingly in the code.\\

\section{Observations on Finite Fields}

\indent Algorithm 6.1 was implemented in a computer algebra system [9] (actually, it was first written in a computer algebra system and then turned into pseudo-code, but whatever). Some results:\\\\
\textbf{Observation 7.1:} In addition to the Parker fields given in the previous section, the only Parker fields of prime order, $p$, with $p<1000$ are $\mathbb{F}_{31},\mathbb{F}_{43},\mathbb{F}_{47},$ and $\mathbb{F}_{67}$.\\\\
\textbf{Corollary 7.1:} The smallest non-Parker field with order $p\equiv3\ (\mod\ 4)$ is $\mathbb{F}_{59}$ with the explicit construction
\begin{equation*}
\begin{tabular}{c|c|c}
  $20^2$ & $12^2$ & $7^2$ \\      \hline
  $2^2$ & $1^2$ & $23^2$ \\      \hline
  $22^2$ & $25^2$ & $29^2$
\end{tabular}\ .
\end{equation*}
\textbf{Conjecture 7.1:} The only Parker fields of prime order are:
$$\mathbb{F}_{2},\ \mathbb{F}_{3},\  \mathbb{F}_{5},\  \mathbb{F}_{7},\  \mathbb{F}_{11},\  \mathbb{F}_{13},\  \mathbb{F}_{17},\  \mathbb{F}_{19},\  \mathbb{F}_{23},\  \mathbb{F}_{31},\ \mathbb{F}_{43},\ \mathbb{F}_{47},\ \text{and}\  \mathbb{F}_{67}.$$
\textbf{Observation 7.2:} The first few finite fields of prime order with more magic squares (up to scaling) than any smaller such finite field are:
\begin{equation*}
\begin{tabular}{c|c}
\text{Field} & \text{\# of magic squares} \\ 
& \text{up to scaling}\\
 \hline
 $\mathbb{F}_{2}$ & $0$ \\
$\mathbb{F}_{29}$ & $2$ \\
$\mathbb{F}_{61}$ & $4$ \\
$\mathbb{F}_{89}$ & $5$ \\
$\mathbb{F}_{97}$ & $6$ \\
$\mathbb{F}_{109}$ & $9$ \\
$\mathbb{F}_{113}$ & $13$ \\
$\mathbb{F}_{137}$ & $18$ \\
$\mathbb{F}_{181}$ & $24$ \\
$\mathbb{F}_{193}$ & $28$ \\
\end{tabular}
\end{equation*}
\textbf{Conjecture 7.2:} The only Parker fields whose order is strictly a prime \textit{power} are $\mathbb{F}_{9},\mathbb{F}_{25},\mathbb{F}_{27},$ and $\mathbb{F}_{243}$.\\

Note that conjectures 7.1 and 7.2 taken together assert that exactly $17$ finite fields are Parker. We could only think to justify the claim with an argument from probability -- which more or less means we couldn't justify the claim. For what it's worth, we point out that the 2 conjectures of this section agree, but are not proven by Theorem 28 in [10] (a finite field is ``Parker" in Labruna's terminology if it contains a magic square of squares ``of order $9$").

\section{A Search Algorithm for Rings}

 In this section, we address which rings of the form $\mathbb{Z}/n\mathbb{Z}$ admit a $3\times3$ magic square of squares -- that is, which such rings are not Parker. \\
\indent Our first objective will be modifying Algorithm 6.1 to work on rings. As we did before, first the algorithm pseudo-code will be given and then an explanation.\\\\
\textbf{Algorithm 8.1:}\\
\# Input: An integer, $n$.\\
\# Output: Set of all tuples $(a^2,b^2,...,i^2)$ forming magic squares over $\mathbb{Z}/n\mathbb{Z}\}$\\
\# $\quad$ up to scaling by $(\mathbb{Z}/n\mathbb{Z})^\times$.\\
\textbf{function} msos\_ring($n$):\\
\indent SQUARES $\leftarrow \{x^2:x^2\in \mathbb{Z}/n\mathbb{Z}\}$\\
\indent MSOS $\leftarrow \{\}$\\
\indent \textbf{for} $m|n$:\\
\indent\indent $e\leftarrow \overline{m}\quad\quad$ \# The overline, $\overline{\ \cdot\ }$, is used to denote a residue in $\mathbb{Z}/n\mathbb{Z}$.\\
\indent\indent SEQUENCES $\leftarrow \{\}$\\
\indent\indent \textbf{for} $\{a^2,i^2\}\subset$ SQUARES:\\
\indent\indent\indent \textbf{if} $a^2+i^2\not=2e^2$: \textbf{continue}\\
\indent\indent\indent \textbf{for} $\{c^2,g^2\}\in$ SEQUENCES:\\
\indent\indent\indent\indent $B\leftarrow 3e^2-a^2-c^2$ \\
\indent\indent\indent\indent $D\leftarrow 3e^2-a^2-g^2$ \\
\indent\indent\indent\indent $F\leftarrow 3e^2-c^2-i^2$ \\
\indent\indent\indent\indent $H\leftarrow 3e^2-g^2-i^2$ \\
\indent\indent\indent\indent \textbf{if} $\{B,D,F,H\}\not\subset$ SQUARES: \textbf{continue}\\
\indent\indent\indent\indent \textbf{if} $|\{a^2,B,c^2,D,e^2,F,g^2,H,i^2\}|\not=9$: \textbf{continue}\\
\indent\indent\indent\indent MSOS $\leftarrow \text{MSOS}\cup \{(a^2,B,c^2,D,e^2,F,g^2,H,i^2)\}$\\
\indent\indent\indent SEQUENCES $\leftarrow\text{SEQUENCES}\cup\{\{a^2,i^2\}\}$\\
\indent \textbf{return} MSOS\\\\
\indent Implementation in the SageMath language is available at [11].\\
\indent For fields, we reduced to the cases in which the central entry of the magic square was either $0$ or $1$. A similar approach is taken for $\mathbb{Z}/n\mathbb{Z}$. Consider the orbits formed under the multiplicative action of $(\mathbb{Z}/n\mathbb{Z})^\times$. There is one orbit corresponding to each divisor of $n$. Thus any magic square of squares over $\mathbb{Z}/n\mathbb{Z}$ may be scaled such that $e$ is the residue of some divisor of $n$.\\

\section{Observations on Rings}

\indent The code from the previous section was again implemented in a computer algebra system [9] resulting in \\\\
\textbf{Observation 9.1:} The smallest ring of the form  $\mathbb{Z}/n\mathbb{Z}$ which admits a magic square of squares is $\mathbb{Z}/27\mathbb{Z}$.\\\\
\textbf{Observation 9.2:} The first few rings of the form $\mathbb{Z}/n\mathbb{Z}$ with more magic squares (up to scaling by $(\mathbb{Z}/n\mathbb{Z})^\times$) than any smaller such ring are:
\begin{equation*}
\begin{tabular}{c|c}
\text{Ring} & \text{\# of magic squares} \\ 
& \text{up to scaling}\\
 \hline
$\mathbb{Z}/27\mathbb{Z})$ & $3$ \\
$\mathbb{Z}/29\mathbb{Z})$ & $7$ \\
$\mathbb{Z}/37\mathbb{Z})$ & $9$ \\
$\mathbb{Z}/53\mathbb{Z})$ & $13$ \\
$\mathbb{Z}/54\mathbb{Z})$ & $36$ \\
$\mathbb{Z}/58\mathbb{Z})$ & $56$ \\
$\mathbb{Z}/74\mathbb{Z})$ & $72$ \\
$\mathbb{Z}/101\mathbb{Z})$ & $75$ \\
$\mathbb{Z}/106\mathbb{Z})$ & $104$ \\
$\mathbb{Z}/122\mathbb{Z})$ & $240$ \\
$\mathbb{Z}/162\mathbb{Z})$ & $576$ \\
$\mathbb{Z}/202\mathbb{Z})$ & $604$ \\
\end{tabular}\quad\quad
\begin{tabular}{c|c}
\text{Ring} & \text{\# of magic squares} \\ 
& \text{up to scaling}\\
 \hline
$\mathbb{Z}/218\mathbb{Z})$ & $680$ \\
$\mathbb{Z}/226\mathbb{Z})$ & $832$ \\
$\mathbb{Z}/274\mathbb{Z})$ & $1296$ \\
$\mathbb{Z}/314\mathbb{Z})$ & $1304$ \\
$\mathbb{Z}/346\mathbb{Z})$ & $2112$ \\
$\mathbb{Z}/362\mathbb{Z})$ & $2262$ \\
$\mathbb{Z}/386\mathbb{Z})$ & $2260$ \\
$\mathbb{Z}/394\mathbb{Z})$ & $2420$ \\
$\mathbb{Z}/458\mathbb{Z})$ & $2904$ \\
$\mathbb{Z}/466\mathbb{Z})$ & $2972$ \\
$\mathbb{Z}/482\mathbb{Z})$ & $3860$ \\
$\mathbb{Z}/486\mathbb{Z})$ & $6120$ \\
\end{tabular}
\end{equation*}
\textbf{Conjecture 9.2:} If  $\mathbb{Z}/n\mathbb{Z}$ has more magic squares (up to scaling) than any smaller such ring, then $n=2p$ for some prime, $p\equiv1(\mod\ 4)$, with exactly 7 exceptions: $n=27,29,37,53,54,101,$ and $162$.\\\\
\textbf{Observation 9.3:} The only Parker $\mathbb{Z}/n\mathbb{Z}$ for odd $n$ in the range $100<n<1000$ are $n=129,141,147,$ and $201$.\\\\
\textbf{Conjecture 9.3:} $\mathbb{Z}/201\mathbb{Z}$ is the largest Parker ring of the form $\mathbb{Z}/n\mathbb{Z}$ where $n$ is odd.\\\\
\textbf{Observation 9.4:} The only Parker $\mathbb{Z}/n\mathbb{Z}$ for $n$ divisible by $4$ in the range $1000<n<3000$ are $n=1032, 1072, 1104, 1128, 1488,1608,2064,$ and $2256$. Interestingly, each of these integers is of the form $2^a3^bp$ with a prime $p\equiv 3\ (\text{mod}\ 4)$.
$$1032=2^3\cdot 3\cdot 43,\quad 1072=2^4\cdot 67,\quad 1104=2^4\cdot 3\cdot 23,\quad 1128=2^3\cdot 3\cdot 47,$$
$$1488=2^4\cdot 3\cdot 31,\quad 1608=2^3\cdot 3\cdot 67,\quad 2064=2^4\cdot 3\cdot 43,\quad\text{and}\quad 2256=2^4\cdot 3\cdot 47.$$
Further searches targeting integers of this form revealed only that $\mathbb{Z}/3216\mathbb{Z}$ is Parker (which indeed fits the form: $3216=2^4\cdot 3\cdot 67$). It is surprising though that such large Parker rings exist.\\\\
\textbf{Conjecture 9.4:} $\mathbb{Z}/3216\mathbb{Z}$ is the largest Parker ring of the form $\mathbb{Z}/n\mathbb{Z}$.\\\\
\indent Note a proof of the existence of infinitely many Parker rings of the form $\mathbb{Z}/n\mathbb{Z}$ would prove that no magic square of distinct squares exists over the integers. If, instead, such a square existed and $N$ were its largest entry, then $\mathbb{Z}/n\mathbb{Z}$ would be non-Parker for $n>N$. And therefore only finitely many such rings could be Parker.

$$\text{REFERENCES}$$
[1] Christian Boyer, \textit{Magic Square of Squares},\\
\indent http://www.multimagie.com/English/SquaresOfSquares.htm
\newline
[2] Matt Parker, \textit{The Parker Square}, Numberphile,\\
\indent $[$interview by Brady Haran$],$ \\
\indent https://www.youtube.com/watch?v=aOT\_bG-vWyg
\newline
[3] John P. Robertson, \textit{Magic Squares of Squares},\\ 
\indent Mathematics Magazine, vol. 69, no. 4, 1996, pp. 289–293.\\ 
\indent JSTOR, www.jstor.org/stable/2690537.
\newline
[4] MathWorld, \textit{Congruum Problem},\\
\indent http://mathworld.wolfram.com/CongruumProblem.html
\newline
[5] Eknath Ghate, \textit{The Kronecker-Weber Theorem}, \\ 
\indent Summer School on Cyclotomic fields, Pune, June 7-30, 1999
\newline
[6] Keith Conrad, \textit{The Gaussian Integers},\\
\indent https://kconrad.math.uconn.edu/blurbs/
\newline
[7] Artin, \textit{Algebra}
\newline
[8] Math Pages, \textit{Magic Square of Squares},\\
\indent https://www.mathpages.com/home/kmath417/kmath417.htm
\newline
[9] SageMath, the Sage Mathematics Software System (Version 8.4),\\
\indent The Sage Developers, 2015, http://www.sagemath.org.
\newline
[10] Giancarlo Labruna, \textit{Magic Squares of Squares of Order Three Over}\\
\indent\textit{Finite Fields}, https://digitalcommons.montclair.edu/etd/138/
\newline
[11] Code repository: https://github.com/onnomc/parker-ring-search

\end{document}